\documentclass[11pt]{amsart}
\theoremstyle{plain}

\begin{document}
\large

\title[$L^2$ Forms and Ricci flow ]
{$L^2$ p-Forms and Ricci flow with bounded curvature on manifolds }
\author{Li MA and Baiyu Liu}

\address{Department of mathematical sciences \\
Tsinghua university \\
Beijing 100084 \\
China}
\thanks{The research is partially supported by the National Natural Science
Foundation of China 10631020 and SRFDP 20060003002}
\email{lma@math.tsinghua.edu.cn} \dedicatory{}
\date{Jan. 5th, 2007}

\keywords{Ricci flow, forms, monotonicity.} \subjclass{53C.}

\begin{abstract}
In this paper, we study the evolution of $L^2$ p-forms under Ricci
flow  with bounded curvature on a complete non-compact or a compact
Riemannian manifold. We show that under the curvature operator bound
condition on such a manifold, the weighted $L^2$ norm of a smooth
p-form is non-increasing along the Ricci flow. The weighted
$L^{\infty}$ norm is showed to have monotonicity property too.
\end{abstract}
\maketitle

\section{Introduction}
In this paper we study the evolution of a p-form under the Ricci
flow introduced by R.Hamilton in 1982 (\cite{H95}). To understand
the change of the DeRham cohomology of the manifold under Ricci
flow, we need to compute the heat equation for p-forms. Then we try
to use some trick from the paper \cite{L} to get some monotonicity
results.

By definition, the Ricci-Hamilton flow on a manifold $M$ of
dimension $n$ is
  the evolution equation for Riemannian metrics:
  $$
\partial_tg_{ij}=-2R_{ij},\;\;\;\mbox{on $M_T:=M\times [0,T)$}
  $$
where $R_{ij}$ is the Ricci tensor of the metric $g:=g(t)=(g_{ij})$
in local coordinates $(x^i)$ and $T$ is the maximal existing time
for the flow. Given an initial complete Riemannian metric of bounded
curvature, the existence of Ricci flow with bounded sectional
curvature on a complete non-compact Riemannian manifold had been
established by Shi \cite{Shi89} in 1989. This is a very useful
result in
 Riemannian Geometry. Interestingly, the maximum principle of heat
 equation is true on such a flow, see \cite{Shi89}. Then we can easily show that
 the Ricci flow
preserves the property of nonnegative scalar curvature (see also
\cite{H95}). Given a smooth $L^2$ p-form $\phi$ with compact support
on a Riemannian manifold $(M,g)$. Recall that its $L^2$ norm is
defined by
$$
||\phi||_{L^2g}=(\int_M|\phi|^2_g(x)dv_g)^{1/2},
$$
and the $L^{\infty}$ norm is defined as $$
||\phi||_{L^{\infty}g}=\sup_{x\in M}|\phi(x)|_{g(x)}
$$
Assume that $d\phi=0$. Let $\Phi=[\phi]$ be the $L^2$ cohomology
class of the form $\phi$ in $(M,g)$. Define
$$
||\Phi||_{L^2g}=\inf_{\varphi\in \Phi}||\varphi||_{L^2g}
$$
and
$$||\Phi||_2(t)=||\Phi||_{L^2g(t)}$$
for the flow $\{g(t)\}$.  It is well-known that $||\Phi||_{L^2g}$ is
a norm on $H^1_{dR}(M,\mathbf{R})$. We denote by $d_g(x,y)$ the
distance of two points $x$ and $y$ in $(M,g)$.

Our new results are the following
\newtheorem{Theorem}{Theorem}
\begin{Theorem}\label{Theorem:thm1} Let $(M, g_{0})$ be a compact or complete noncompact Riemannian
manifold with non-negative scalar curvature. Assume that $g(t)$ is a
Ricci flow with bounded curvature on $[0,T)$ with initial metric
$g(0)=g_{0}$ on M. For $t \in [0,T)$ the Ricci curvature
satisfies \\
\begin{displaymath}
R_{ij}\eta ^{i}\eta ^{j}\geq 0
\end{displaymath}
and it holds the curvature pinching condition
\begin{equation}\label{eq8} W(t)+\frac {2R(t)}{(n-1)(n-2)}
\leq \frac {4}{n-2}L(t),
\end{equation}
where $R(t)$ is the scalar curvature of the flow $(g(t))$, and
$$
W(t)=\sup_{\xi } \frac {\lvert W_{ijkl}(t) \xi ^{ij}\xi ^{kl}
\lvert}{\xi ^{ij} \xi _{ij}},$$
$$(\xi ^{ij}=-\xi ^{ji})$$ with
$$W_{ijkl}=R_{ijkl}-\frac
{1}{n-2}(R_{ik}g_{jl}-R_{il}g_{jk}+R_{jl}g_{ik}-R_{jk}g_{il})+\frac
{R}{(n-1)(n-2)}(g_{ik}g_{jl}-g_{il}g_{jk})$$ is the Weyl conformal
curvature tensor. $L(t)$ is the smallest eigenvalue of the matric
$R_{ij}(t)$.\\

 Then for a $L^{2}$ p-form $\xi$ , we have
\begin{displaymath}
\lvert\lvert \xi \rvert \rvert_{L^{2}g(t)} \leq \lvert\lvert \xi \rvert
\rvert_{L^{2}g(s)}, for \ t>s ,
\end{displaymath}
alone the Ricci flow $g(t)$. Similarly, we have the $L^{\infty}$
monotonicity of the p-form heat flow along the Ricci flow.
\end{Theorem}

It is easy to see that the $L^2$ monotonicity gives the monotonicity
of De Rham cohomology class of a closed p-form with compact support.
So we shall not state the corresponding result for De Rham
cohomology class. We remark that the pinching curvature condition in
the theorem above is not nature since it may not be preserved along
the Ricci flow. However, this is a classical condition used in the
book of Bochner and Yano \cite{B} ( see page 89 of the \cite{B}),
which is quite
similar to ours. \\

As a comparison, we would like to mention a pinching result of
G.Huisken \cite{H85}. It is well known that the curvature tensor
$Rm=\{ R_{ijkl} \}$ of a Riemannian manifold can be decomposed into
three orthogonal components which have the same symmetries as $Rm$:
$$Rm=W+V+U .$$
Here $W=\{ W_{ijkl} \}$ is the Weyl conformal curvature tensor,
whereas $V=\{ V_{ijkl} \}$ and $U=\{ U_{ijkl} \}$ denote the
traceless Ricci part and the scalar curvature part respectively. The
following pointwise pinching condition was proposed by Huisken in
\cite{H85} (see also the works of C.Margerin \cite{Ma86} and
S.Nishikawa\cite{Ni86}):
\begin{equation}\label{pinch} \lvert W \lvert
^{2}+\lvert V \lvert ^{2} < \delta _{n}  \lvert U \lvert ^{2} ,
\end{equation}
with
\begin{equation}\label{deltan}\delta _{4}=\frac {1}{5},\qquad \delta _{5}=\frac {1}{10},\qquad \delta _{n}=\frac {2}{(n-2)(n+1)} ,\quad n\geq 6,
\end{equation}
and define the norm of a tensor:
$$\lvert T\lvert^{2}=\lvert T_{ijkl}\lvert^{2}=g^{im}g^{jn}g^{kp}g^{lq}T_{ijkl}T_{mnpq}=T_{ijkl}T^{ijkl}. $$

Then we have the following result of G.Huisken:
\newtheorem{Theorem:compact}{Theorem}
\begin{Theorem}\label{Theorem:compact}
Let $n\geq 4$. Suppose $(M^{n}, g_{0})$ is an n-dimensional smooth
compact Riemannian manifold with positive and bounded scalar
curvature and satisfies the pointwise pinching condition
(\ref{pinch}).
\\
Then $M^n$ is diffeomorphic to the sphere $S^n$ or a quotient space
of $S^n$ by a group of fixed point free isometries in the standard
metric.
\end{Theorem}
Note that the pinching curvature conditions (even the assumption
about behavior at infinity see \cite{D94} and \cite{DM05}) in
Theorem 2 and Theorem 3 are preserved along the Ricci flow, see
\cite{H85} \cite{Ma86} and \cite{Ni86}. So our pinching condition in
Theorems 2 is more nature. We point out that the pinching condition
(\ref{pinch}) gives positive curvature operator. Recently, the deep
work of C.Boehm and B.Wilking \cite{BW06} proved the same result for
positive curvature operator case.

In the following we just try to give another way to understand the
monotonicity of the norms of closed forms under Ricci flow for
general lower bound curvature operator case.

\newtheorem{Theorem:3compact}{Theorem}
\begin{Theorem}\label{Theorem:3compact}
Let $n\geq 4$. Suppose $(M^{n}, g_{0})$ is an n-dimensional smooth
compact Riemannian manifold.
 Assume that $g(t)$ is a Ricci flow with its curvature operator bounded from below
 by the constant $2k$ on
$[0,T)$ with
initial metric $g(0)=g_{0}$ on M.\\

 Then for a $L^{2}$ p-form $\xi$ on $(M, g_0)$ , we have
\begin{displaymath}\lvert\lvert e^{kp(p-1)t}\xi \rvert \rvert_{L^{2}g(t)}
\leq \lvert\lvert e^{kp(p-1)s}\xi \rvert \rvert_{L^{2}g(s)}, for \
t>s ,
\end{displaymath}
alone the Ricci flow $g(t)$. Similarly, we have the $L^{\infty}$
monotonicity of the p-form heat flow along the Ricci flow.
\end{Theorem}

Note that R.Hamilton \cite{H95} proved that the non-negative
curvature operator condition is preserved along the Ricci flow on
compact Riemannian manifolds. With the help of the curvature decay
estimate, one can show that similar result to Theorem 3 is also true
in complete non-compact Riemannian manifolds. However, we omit the
detail of the proof here since the argument is similar to Theorem 3.
\newtheorem{Theorem:non-compact}{Theorem}
\begin{Theorem}\label{Theorem:non-compact}
Let $n\geq 4$. Suppose $(M^{n}, g_{0})$ is an n-dimensional smooth
complete noncompact Riemannian manifold with bounded curvature
operator and with curvature decay condition as in \cite{DM05}.
Assume that $g(t)$ is a Ricci flow with its curvature operator
bounded from below
 by the constant $2k$ on $[0,T)$
with
initial metric $g(0)=g_{0}$ on M. \\

 Then for a $L^{2}$ p-form $\xi$ on $(M, g_0)$, we have
\begin{displaymath}\lvert\lvert e^{kp(p-1)t}\xi \rvert \rvert_{L^{2}g(t)}
\leq \lvert\lvert e^{kp(p-1)s}\xi \rvert \rvert_{L^{2}g(s)}, for \
t>s ,
\end{displaymath}
alone the Ricci flow $g(t)$. Similarly, we have the $L^{\infty}$
monotonicity of the p-form heat flow along the Ricci flow.
\end{Theorem}

\section{Basic formulae from Riemannian Geometry}

Some basic materials in Riemannian geometry are stated here, to the
extent that will serve as computational notations the later
sections. Readers who are interested in pursuing further along the
line are referred to the book by Yano and Bochner \cite{B} and the
paper by Huisken \cite{H85}. However, we make a caution that we use
modern convention from the book of \cite{C}.

Consider an n-dimensional Riemannian manifold $M^{n}$ with the metric $(g_{ij})$.
Denote by $(g^{ij})=(g_{ij})^{-1}$ and $\Gamma ^{i} _{jk}$ the Christoffel symbols.\\

For a scalar $f(x)$, the covariant derivative of $f(x)$ is given by
\begin{displaymath}f_{;j}=\frac{\partial f}{\partial x^{j}}
\end{displaymath}
and the second covariant derivative is given by
\begin{displaymath}f_{;j;k}=\frac{\partial^{2} f}{\partial x^{j} \partial
x^{k}}-\frac{\partial f}{\partial x^{i}}\Gamma ^{i} _{jk}.
\end{displaymath}

Thus, we see that $f_{;j;k}=f_{;k;j}$. \\

However, for vectors and tensors, successive covariant
differentiations are not commutative in general. For example, for a
contravariant vector $v^{i}$, we obtain
\begin{equation}\label{eq1} v^{i}_{;k;l}-v^{i}_{;l;k}=-v^{j}
R^{i}_{klj} \ ,
\end{equation}
where
\begin{equation} \label{eq4} R^{l}_{ijk}=\frac{\partial \Gamma ^{l}_{jk}}{\partial
x^{i}}-\frac{\partial \Gamma ^{l}_{ik}}{\partial x^{j}}+\Gamma
^{l}_{im}\Gamma ^{m}_{jk}-\Gamma ^{l}_{jm}\Gamma ^{m}_{ik}.
\end{equation}

Similarly, for a covariant vector $v_{j}$, then we have
\begin{equation}\label{eq2}
v_{k;i;j}-v_{k;j;i}=v_{l}R^{l}_{ijk} \ ,
\end{equation}
and if we take a general tensor $T^{i}_{j k}$ for example, then we
have

\begin{equation} \label{eq3}
T^{i}_{j k;l;m}-T^{i}_{j k;m;l}=-T^{s}_{j k} R^{i}_{lms}+T^{i}_{s k}
R^{s}_{lmj}+T^{i}_{j s} R^{s}_{lmk} .
\end{equation}

Formulas (\ref{eq1}),(\ref{eq2}) and (\ref{eq3}) are called the
Ricci formulas. \\

From the curvature tensor $R^{i}_{jkl}$, we get, by contraction,
\begin{displaymath}R_{jk}=R^{i}_{ijk} ,
\end{displaymath}
moreover, from $R_{jk}$, by multiplication by $g^{jk}$ and by
contraction, we get
\begin{displaymath}R=g^{jk}R_{jk} .
\end{displaymath}

$R_{jk}$ and $R$ are called  Ricci tensor and curvature
scalar of the metric $g$ respectively. \\

From the definition (\ref{eq4}) of $R^{i}_{jkl}$ , it is easily seen
that $R^{i}_{jkl}$ satisfies the following algebraic identities:
\begin{displaymath}R^{i}_{jkl}=-R^{i}_{jlk} ,
\end{displaymath}
and
\begin{equation} \label{eq5} R^{i}_{jkl}+R^{i}_{klj}+R^{i}_{ljk}=0 .
\end{equation}

Consequently, from
(\ref{eq5}), we obtain $R_{jk}=R_{kj}$.\\

If we put
\begin{displaymath}R_{ijkl}=g_{hl} R^{h}_{ijk} ,
\end{displaymath}
then $R_{ijkl}$ satisfies
\begin{displaymath}R_{ijkl}=-R_{jikl} ,
\end{displaymath}
and
\begin{equation} \label{eq6} R_{ijkl}+R_{jkil}+R_{kijl}=0 .
\end{equation}

Equations (\ref{eq5}) and (\ref{eq6}) are called the first Bianchi
identities. \\

Moreover, applying the Ricci formula and calculating the covariant
compontents $R_{ijkl}$, we get
\begin{displaymath}R_{ijkl}=-R_{ijlk} ,
\end{displaymath}

and
\begin{displaymath}R_{ijkl}=R_{klij} .
\end{displaymath}\\

It is also to be noted that
\begin{equation}\label{bianch2} R^{i}_{jkl;m}+R^{i}_{jlm;k}+R^{i}_{jmk;l}=0 ,
\end{equation}

which is called the second Bianchi identity. From (\ref{bianch2}),
we get
\begin{displaymath}2R^{s}_{l;s}=R_{;l} ,
\end{displaymath}
in which $R^{s}_{l}=g^{is}R_{il} .$

Denote by $Rc=\{R_{ij}\}$  and $R$ the Ricci tensor and scalar
curvature. We can write the traceless Ricci part $V=\{V_{ijkl}\}$
and the scalar curvature part $U=\{U_{ijkl}\}$ as follows ( see also
\cite{H85}):

$$U_{ijkl}=\frac{1}{n(n-1)}R(g_{ik}g_{jl}-g_{il}g_{jk}) ,$$
$$V_{ijkl}=\frac{1}{n-2}(\overset {\circ}{R}_{ik}g_{jl}-
\overset {\circ}{R}_{il}g_{jk}-\overset
{\circ}{R}_{jk}g_{il}+\overset {\circ}{R}_{jl}g_{ik}),$$ where
$$\overset {\circ}{R}_{ij}=R_{ij}-\frac{1}{n}Rg_{ij}\ .$$
If we let $$\overset {\circ}{Rm}=\{\overset
{\circ}{R}_{ijkl}\}=\{R_{ijkl}-U_{ijkl}\}=\{V_{ijkl}+W_{ijkl}\} ,$$
then
$$\vert \overset {\circ}{Rm}\vert^{2}=\vert W \vert^{2}+\vert V \vert^{2} ,$$
$$\vert U \vert^{2}=\frac {2}{n(n-1)}R^{2} ,$$
$$\vert Rm \vert^{2}=\vert \overset {\circ}{Rm}\vert^{2}+\vert U \vert^{2} .$$

\section{Proof of Theorem ~\ref{Theorem:thm1}}

We first set up a key lemma which is useful in the proof of all the
Theorems above.
\newtheorem{lemma}{Lemma}
\begin{lemma}Let $\xi _{i_{1} \dots i_{p}}(x,t)$ be an
anti-symmetric covariant vector for all the time t, and satisfying
the heat equation
\begin{displaymath}  \frac{\partial \xi _{i_{1} \dots i_{p}}}{\partial
t}=\Delta _{d} \, \xi _{i_{1} \dots i_{p}} .
\end{displaymath}

Then,
\begin{equation}\label{eq7} \frac{\partial}{\partial t} \lvert \xi
\rvert ^{2}=\Delta \lvert \xi \rvert ^{2} -2\xi ^{i_{1} \dots
i_{p};j} \xi _{i_{1} \dots i_{p};j}-p(p-1)\xi _{iji_{1} \dots
i_{p-2}}R^{ij}_{\phantom{ij}kl} \xi ^{kli_{1}\dots i_{p-2}} .
\end{equation}

Here $\Delta _{d} =\delta  d+d \delta$ is the Hodge-DeRham Laplacian
of $g(t)$, $\Delta$ is the rough Laplacian in the sense of \cite{C},
and $R^{ij}_{\phantom{ij}kl}=g^{ip}g^{jq}R_{pqkl}$.

\end{lemma}

Proof: Recall that for a covariant vector $\xi _{i_{1} \dots
i_{p}}$, we have

\begin{displaymath}\Delta
_{d} \, \xi _{i_{1} \dots i_{p}}=g^{jk} \xi _{i_{1} \dots
i_{p};j;k}- \sum^{1\dots p}_{s}\xi _{i_{1} \dots i_{s-1} a i_{s+1}
\dots i_{p}} R^{a}_{\phantom{a} i_{s}}- \sum ^{1 \dots p}_{s<t} \xi
_{i_{1} \dots i_{s-1} a i_{s+1} \dots i_{t-1} b i_{t+1} \dots i_{p}}
R^{ab}_{\phantom{ab} i_{s} i_{t}} ,
\end{displaymath}
(see page 74 in \cite{B},
 and please note that we have used the same notations as in \cite{B}). \\

Along the Ricci-Hamilton flow, we have
$$\frac{\partial g^{ij}}{\partial t}=2R^{ij}.$$

Then,
\begin{eqnarray}\nonumber \frac{\partial}{\partial t} \lvert \xi \rvert ^{2} & = &
\frac{\partial}{\partial t} (g^{i_{1}k_{1}} \dots g^{i_{p}k_{p}}\xi
_{i_{1} \dots i_{p}}\xi _{k_{1} \dots k_{p}})\\\nonumber & = &
\sum^{1 \dots p}_{a} g^{i_{1}k_{1}} \dots g^{i_{a-1}k_{a-1}}
2R^{i_{a}k_{a}} g^{i_{a+1}k_{a+1}} \dots g^{i_{p}k_{p}}\xi _{i_{1}
\dots i_{p}}\xi _{k_{1} \dots k_{p}}
\\ \nonumber
&  &  + 2g^{i_{1}k_{1}} \dots g^{i_{p}k_{p}} (\Delta _{d} \,
\xi _{k_{1} \dots k_{p}}) \xi _{i_{1} \dots i_{p}} \\
\nonumber & = & 2p R^{i_{1}k_{1}} g^{i_{2}k_{2}}\dots
g^{i_{p}k_{p}}\xi _{i_{1} \dots i_{p}}\xi _{k_{1} \dots k_{p}}
+2g^{i_{1}k_{1}} \dots g^{i_{p}k_{p}} (\Delta \xi _{k_{1} \dots
k_{p}}) \xi _{i_{1} \dots i_{p}} \\\nonumber &  & -2\sum^{1\dots
p}_{s}\xi ^{k_{1} \dots k_{p}} \xi_{k_{1} \dots
k_{s-1} a k_{s+1}\dots k_{p}}R^{a}_{\phantom{a} k_{s}}\\
\nonumber &  & -2\sum ^{1 \dots p}_{s<t} \xi _{k_{1} \dots k_{s-1} a
k_{s+1} \dots k_{t-1} b k_{t+1} \dots k_{p}} R^{ab}_{\phantom{ab}
k_{s} k_{t}}\xi ^{k_{1} \dots k_{p}}\\
\label{eq9}& = & 2p R_{ij} \xi ^{j i_{2} \dots i_{p}}\xi ^{i}
_{\phantom{i} i_{2} \dots i_{p}}+2 \xi ^{k_{1} \dots
k_{p}}(\bigtriangleup \xi _{k_{1} \dots k_{p}})-2p \xi ^{ji_{2}
\dots i_{p-1}} \xi _{ii_{2} \dots i_{p-1}}R^{i}_{\phantom{i}j}\\
\nonumber & & -p(p-1)\xi _{iji_{1} \dots
i_{p-2}}R^{ij}_{\phantom{ij}kl} \xi ^{kli_{1}\dots i_{p-2}} .
\end{eqnarray}
By calculations, we obtain
\begin{equation}\label{laplacian} \Delta \lvert \xi \rvert ^{2}=2(\xi^{i_{1} \dots
i_{p}} \Delta \xi_{i_{1} \dots i_{p}}+\xi^{i_{1} \dots i_{p};j}
\xi_{i_{1} \dots i_{p};j}),
\end{equation}
for which $\Delta f=g^{ij}f_{;i;j}$ , (f is a scalar field) . \\

 Putting (\ref{laplacian}) into (\ref{eq9}), we get (\ref{eq7})
 $$\frac{\partial}{\partial
t} \lvert \xi \rvert ^{2}=\Delta \lvert \xi \rvert ^{2} -2\xi
^{i_{1} \dots i_{p};j} \xi _{i_{1} \dots i_{p};j}-p(p-1)\xi
_{iji_{1} \dots i_{p-2}}R^{ij}_{\phantom{ij}kl} \xi ^{kli_{1}\dots
i_{p-2}}.$$

We are done since this is what we wanted in the lemma.

We now give a proof of Theorem 1.

 {\bf Proof of Theorem
~\ref{Theorem:thm1}}: Using the curvature decomposition we have the
following:
\begin{eqnarray}\nonumber & \quad{-} &\xi _{iji_{1} \dots
i_{p-2}}R^{ij}_{\phantom{ij}kl} \xi ^{kli_{1}\dots i_{p-2}}
\\ \nonumber& = & - \xi ^{iji_{1} \dots i_{p-2}}R_{ijkl} \xi
^{kl}_{\phantom{kl} i_{1}\dots
i_{p-2}}\\
\nonumber & = & \xi ^{iji_{1} \dots i_{p-2}} \xi ^{kl}_{\phantom{kl}
i_{1}\dots i_{p-2}} ( -W_{ijkl} +\frac
{R(t)}{(n-1)(n-2)}(g_{ik}g_{jl}-g_{il}g_{jk})\\
\nonumber &  & -\frac
{1}{n-2}(R_{ik}g_{jl}-R_{il}g_{jk}+R_{jl}g_{ik}-R_{jk}g_{il})) \\
\nonumber & = & -W_{ijkl}\xi ^{iji_{1} \dots i_{p-2}} \xi
^{kl}_{\phantom{kl} i_{1}\dots i_{p-2}}-\frac {4}{n-2} R_{il} \xi
^{iji_{1} \dots i_{p-2}} \xi ^{l}_{\phantom{l} ji_{1}\dots
i_{p-2}}\\
\nonumber &  & +\frac {2R(t)}{(n-1)(n-2)}\xi ^{iji_{1} \dots
i_{p-2}}
\xi _{iji_{1}\dots i_{p-2}}\\
\nonumber & \leq & \left( W(t)-\frac {4}{n-2}L(t)+\frac
{2R(t)}{(n-1)(n-2)} \right) \xi ^{i_{1} \dots i_{p}} \xi
_{i_{1}\dots i_{p}}
\end{eqnarray}
Recall that the assumption (\ref{eq8}) $$W(t)+\frac
{2R(t)}{(n-1)(n-2)}\leq \frac {4}{n-2}L(t),$$
 for (\ref{eq7}) , and we derive
$\frac{\partial}{\partial t} \lvert \xi \rvert ^{2} \leq \Delta
\lvert \xi \rvert ^{2}$. Then, Lemma 4 in \cite{L} implies the $L^2$
monotone result, and the Maximum principle \cite{Shi89} implies the
$L^{\infty}$ monotonicity.
This completes the proof of Theorem 1.\\

\section{Proof of Theorem ~\ref{Theorem:3compact}}

In this section, we plan to prove Theorems 3.

{\bf Proof of Theorem ~\ref{Theorem:3compact}}: First, the existence
of the heat flow of the p-form along the Ricci flow follows almostly
from Gaffney \cite{Ga54}, so we omit the detail. Then we note that
the $L^2$ property of the p-form $\xi$ is preserved along the Ricci
flow. Note that when $k=0$, the positivity of curvature operator is
preserved as long as the solution of
  the evolution equation for the Ricci-Hamilton flow exists. Recall
  that our curvature operator is bounded from below along the Ricci
  flow.
Then we have,
$$\xi _{iji_{1} \dots
i_{p-2}}R^{ij}_{\phantom{ij}kl} \xi ^{kli_{1}\dots i_{p-2}} \geq
2k|\xi|^2, \qquad for \quad t \in [0,T)  \ ,$$ which in turn, by the
equation (\ref{eq7}), gives us that
$$
\frac{\partial}{\partial t} \lvert \xi \rvert ^{2}\leq\Delta \lvert
\xi \rvert ^{2}-2kp(p-1)|\xi|^2.
$$
In the other word, we have
$$
\frac{\partial}{\partial t} \lvert e^{kp(p-1)t}\xi \rvert
^{2}\leq\Delta \lvert e^{kp(p-1)t}\xi \rvert ^{2}.
$$

 Hence, using the same argument as in
(\cite{L}), we have proved the monotonicity of the weighted $L^2$
norm of the p-form. By the Maximum principle \cite{Shi89} we have
the monotonicity of the $L^{\infty}$ norm of the p-form, and then
the result of Theorem 3 has been completely proved.\\


\begin{thebibliography}{99}
\bibitem{B}   S. Bochner \& K. Yano, \emph{Curvature and Bitti Numbers}, Princeton, New Jersey Princeton University
Press, (1953) 16-19, 74, 84.

\bibitem{BW06}
C.Boehm and B.Wilking, \emph{Manifolds with positive curvature
operators are space forms}, math.DG/0606187, 2006.

\bibitem{ZXP2000} Bing-Long Chen \& Xi-Ping Zhu, \emph{ Complete
Riemannian manifolds with pointwise pinched curvature}, Inventiones
Mathematicae,
 Vol.140, Number 2, (2000) 423-452.

\bibitem{C} Bennett Chow \& Dan Knopf, \emph{The Ricci Flow: An
Introduction}, Mathematical Surveys and Monographs, v.110, (2004)
284.

\bibitem{DM05}
X.Z.Dai and L.Ma,  \emph{Mass under the Ricci flow},
math.DG/0510083, 2005, to appear in Comm. Math. Phys., 2007

\bibitem{D94} G.Drees, \emph{Asymptotically flat manifold of non-negative curvature}, Diff. Geom. and its applications.
4(1994)77-90.

\bibitem{Ga54}
M.Gaffney, \emph{The heat equation of Milgram and Rossbloom for open
Riemannian manifolds}, 60(1954)458-466.

\bibitem{H95}
 R.Hamilton, \emph{ The formation of Singularities in the Ricci flow},
 Surveys in Diff. Geom. ,
 Vol.2, (1995)7-136.

\bibitem{H85}
G.Huisken, \emph{Ricci deformation of the metric on a Riemannian
manifold}, J.Diff. Geom. , 21(1985)47-62.

\bibitem{L} Li Ma and Yang Yang, \emph{$L^{2}$ Forms and Ricci flow with
bounded curvature on complete non-compact mainfolds}, Geom Dedicata,
119, (2006)151-158.

\bibitem{Ma86} C. Margerin, \emph{ Pointwise pinched manifolds are spaces
forms}, Proceedings of Symposia in Pure Mathematics, 44,(1986)307-28
    [Arcata: Geometric Measure Theorey and Calculus of Variations, July
1984]

\bibitem{Ni86} S, Nishikawa, \emph{Deformation of Riemannian metrics and
manifolds with bounded curvature ratios},44,(1986)343-52,
    [Arcata: Geometric Measure Theorey and Calculus of Variations, July
1984]

\bibitem{Shi89} W.X.Shi, \emph{ Deforming the metric on complete
Riemannian manifolds}, J. Diff. Geom. , 30, (1989)353-360.



\end{thebibliography}
\end{document}